\documentclass[a4paper,12pt]{amsart}

\usepackage{amsmath}
\usepackage{amssymb}
\usepackage[francais, english]{babel}
\usepackage[T1]{fontenc}
\usepackage{enumerate}

\newcommand{\D}{\mathcal D}
\newcommand{\demo}{\begin{proof}[Démonstration]}
\newcommand{\fdemo}{\end{proof}}
\newcommand{\cc}{\mathbb{C}}
\newcommand{\nn}{\mathbb{N}}
\newcommand{\nm}[1]{\left|#1\right|}

\newcommand{\hh}{\mathbb{H}}
\newcommand{\oh}{0\oplus H}
\newcommand{\ho}{H\oplus 0}

\newcommand{\g}{\mathfrak{g}}
\newcommand{\gk}{\mathfrak{k}}
\newcommand{\gp}{\mathfrak{p}}

\newcommand{\rr}{\mathbb{R}}

\newcommand{\uu}{\mathbb{U}}

\newcommand{\mt}{\mapsto}
\newcommand{\ra}{\rightarrow}

\newcommand{\F}{\mathcal F}

\newcommand{\eh}{\varepsilon}

\newcommand{\ps}[1]{\left\langle #1 \right\rangle}

\newcommand{\ov}[1]{\overline{#1}}

\newcommand{\trans}[1]{{}^\mathrm{t}\!#1}

\newcommand{\ali}{\begin{align}}
\newcommand{\fali}{\end{align}}

\DeclareMathOperator{\id}{id}

\DeclareMathOperator{\rang}{rang}
\DeclareMathOperator{\im}{im}

\DeclareMathOperator{\Det}{Det}

\DeclareMathOperator{\Li}{L}
\DeclareMathOperator{\Sym}{Sym}

\DeclareMathOperator{\rk}{rang}
\DeclareMathOperator{\cd}{codim}

\newtheorem{lem}{Lemme}[section]

\newtheorem{pro}[lem]{Proposition}
\newtheorem{theo}[lem]{Théorème}
\newtheorem{pd}[lem]{Proposition-Définition}
\newtheorem{defi}[lem]{Définition}
\theoremstyle{definition}

\newtheorem{ex}[lem]{Exemple}
\newtheorem{re}[lem]{Remarque}
\newtheorem{pb}[lem]{Problème}
\theoremstyle{remark}
\newtheorem*{nota}{Notations}

\title[L'indice de Maslov]{L'indice de Maslov dans les $JB^*$-triples}

\author{St\'ephane Merigon}
\address{Fachbereich Mathematik, AG AGF\\
Technische Universität Darmstadt\\
Schlossgartenstrasse 7\\
64289 Darmstadt}
\email{merigon@mathematik.tu-darmstadt.de}
\date{\today}
\begin{document}
\selectlanguage{francais}
\begin{abstract}
Soit $E$ un $JB^*$-triple dont l'ensemble des tripotents inversibles 
$\Sigma$ n'est pas vide. Nous construisons un indice invariant par homotopie sur les chemins dans $\Sigma$ qui respectent une condition de type Fredholm par rapport à un tripotent fixé. Cet indice généralise l'indice de Maslov pour la Fredholm-Lagrangienne d'un espace de Hilbert symplectique de dimension infinie défini dans \cite{BOFU}. Lorsque $E$ est de dimension finie, nous relions cet indice à l'indice triple généralisé de \cite{CLOR1, CL} et à l'indice de Souriau généralisé de \cite {CLKO}.
\end{abstract}

\maketitle

\selectlanguage{english}

\vspace{-1cm}
\begin{abstract}
Let $E$ be a $JB^*$-triple whose set of invertible tripotents $\Sigma$ is not empty. We construct a homotopy invariant index for paths in $\Sigma$ that
satisfie a Fredhom type condition with respect to a fixed invertible tripotent. This index generalises the Maslov index for the Fredholm-Lagrangian of an infinite dimensional symplectic Hilbert space defined in \cite{BOFU}. When $E$ is finite dimensional we make the connection with the generalised triple index of \cite{CLOR1, CL} and the generalised Souriau index of \cite {CLKO}.
\end{abstract}

\selectlanguage{francais}

\vspace{0.5cm}

\setcounter{section}{-1}

\section{Introduction}

Dans son traité {\itshape théorie des perturbations et méthodes asymptotiques}, V.P. Maslov introduit un indice pour les chemins dans la Lagrangienne d'un espace symplectique réel de dimension finie qui intervient dans le prolongement de solutions asymptotiques d'équations aux dérivées partielles. Dans \cite{AR1} (voir aussi \cite{AR2}), Arnold clarifie la définition de cet indice.  
Soit $(H,\omega)$ un espace symplectique réel de dimension $2n$ et notons $\Lambda(n)$ sa Lagrangienne.
Pour tout $\lambda\in\Lambda(n)$ et tout $1\leq k\leq n$ posons 
$$\Lambda_\lambda^k(n)=\{\mu\in\Lambda(n)\mid \dim{\mu\cap\lambda}=k\}.$$
Alors
$$\overline{\Lambda_\lambda^1(n)}=\sum_{1\leq k\leq n}{\Lambda_\lambda^k(n)}$$
est un cycle de lieu singulier $\sum_{2\leq k\leq n}{\Lambda_\lambda^k(n)}$.
Il existe sur $H$ un produit scalaire $(.,.)$ et une structure complexe $J$ isométrique tels que 
$$\forall \eta,\xi\in H,\quad \omega(\xi,\eta)=(J\xi,\eta).$$
On peut alors orienter $\overline{\Lambda_\lambda^1(n)}$ transversalement grâce au champ 
$$v(\mu)={\frac{\mathrm{d}}{\mathrm{d}\theta}}_{\vert{\theta=0}}e^{J\theta}\mu,$$
le coté positif étant celui vers lequel $v(\mu)$ est dirigé.
L'indice (par rapport à $\lambda$) d'un chemin $\gamma$ dont les extrémités ne sont pas dans le cycle est par définition l'indice d'intersection de $\gamma$ avec ce cycle : si l'ensemble des points d'intersection de $\gamma$ avec le cycle est fini et contenu dans $\Lambda_\lambda^1(n)$ et si en chacun de ces points $\gamma$ est continûment différentiable alors l'indice de Maslov est le nombre de points où $\gamma$ traverse le cycle dans le sens positif moins le nombre de points où $\gamma$ traverse le cycle dans le sens négatif. Lorsque l'on se restreint aux chemins fermés on obtient un élément du groupe de cohomologie entière $H^1(\Lambda(n),\mathbb{Z})$ qui ne dépend pas de $\lambda$. 

Remarquons que l'orthogonal (pour le produit scalaire) d'un lagrangien $\lambda$ est $\lambda^\bot=J\lambda$. Fixons une base orthogonale de $\lambda$ et identifions $H$ à $\cc^n$ muni de la forme hermitienne $\left\langle.,.\right\rangle=(.,.)-i\omega(.,.)$ par :
\begin{align*}
H\simeq\lambda\oplus\lambda^\perp&\simeq\cc^n\\
\eta\oplus J\xi&\mapsto\eta+i\xi.
\end{align*}
Alors le groupe $U(n)$ des matrices complexes unitaires de taille $n$ agit transitivement sur $\Lambda(n)$ et le stabilisateur de $\lambda$ s'identifie au sous-groupe des matrices réelles $O(n)$ :
$$\Lambda(n)\simeq U(n)/O(n).$$
On peut donc définir une application 
$\Det^2:\Lambda(n)\rightarrow S^1$
et Arnold montre qu'elle induit un isomorphisme des groupes fondamentaux :
$$\pi_1(\Lambda(n))\simeq \pi_1(S^1).$$
Ainsi on a 
$$H_1(\Lambda(n),\mathbb{Z})\simeq \pi_1(\Lambda(n))$$ et il revient donc au même de se donner un générateur de $H^1(\Lambda(n),\mathbb{Z})$ ou un isomorphisme $\pi_1(\Lambda(n))\simeq\mathbb Z$. Arnold montre que l'indice de Maslov coïncide avec l'image réciproque par $\Det^2$ du générateur standard de $\pi_1(S^1)$ (le nombre de tours sur $S^1$ orienté dans le sens trigonométrique).

Motivé par une justification rigoureuse de la méthode de Maslov, Leray  donne une variante de la définition d'Arnold-Maslov (cf. \cite{LE}). L'indice apparaît comme une fonction sur le double produit du revêtement universel de la Lagrangienne et réalise une primitive d'un cocycle défini sur les triplets de lagrangiens appelé indice d'inertie. Enfin Souriau, grâce à une construction explicite du revêtement universel, donne une formule explicite pour la fonction de Maslov (cf. \cite{SO}). 

Dans \cite{BOFU} Booss-Bavnbek et Furutani généralisent l'indice de Maslov pour la Lagrangienne d'un espace de Hilbert symplectique de dimension infinie $H$. Soit $\lambda$ un lagrangien de $H$. L'indice est défini pour les chemins dans la Fredholm-Lagrangienne $\mathcal{F}\Lambda_\lambda$, c'est-à-dire l'ensemble des lagrangiens $\mu$ tels que $(\lambda,\mu)$ est une paire de Fredholm: 
$$\dim \lambda\cap\mu<\infty\quad\text{et}\quad\dim H/(\lambda+\mu)<\infty,$$
et il réalise un isomorphisme entre $\pi_1(\mathcal{F}\Lambda_\lambda)$ et $\mathbb{Z}$.

Dans une autre direction, Jean-Louis Clerc et Bent {\O}rsted ont montré (cf. \cite{CLOR1, CLOR2, CL}) que l'indice triple se généralise naturellement à la frontière de Shilov $S$ d'un domaine borné symétrique de type tube $\mathcal D$, et qu'il permet de caractériser les orbites de triplets transverses de $S$ sous l'action du groupe des automorphismes holomorphes de $\mathcal D$. Puis Clerc et Koufany (cf. \cite{CLKO}) ont construit de deux manières différentes une primitive de l'indice triple sur le revêtement universel de la frontière de Shilov, l'une généralisant la méthode de Souriau et l'autre celle d'Arnold-Maslov.
A la fin des années 70, Kaup et Upmeier ont développé la théorie des domaines borné symétriques dans les espaces de Banach, le résultat principal étant que la catégorie des domaines bornés symétriques est équivalente à celle des $JB^*$-triples. Dans cet article nous construisons l'indice de Maslov pour l'ensemble des tripotents inversibles d'un $JB^*$-triple, en adaptant la construction de Booss-Bavnbek et Furutani.

Le paragraphe 2 présente la structure de $JB^*$-triple et son lien avec les domaines bornés symétriques. Dans le paragraphe 3, nous détaillons l'identification entre la Lagrangienne d'un espace de Hilbert symplectique réel $H_0\oplus H_0$ et l'ensemble des tripotents inversibles du $JB^*$-triple $Sym(H_0\oplus iH_0)$. Dans le paragraphe 4 nous introduisons la définition d'une paire de Fredholm pour deux unités d'un $JB^*$-triple, et l'indice de transversalité d'une telle paire $(x,e)$ et nous étudions comment évolue cet indice lorsque l'on perturbe $x$. Cette étude nous permet de construire dans le paragraphe 5 l'indice de Maslov d'un chemin $t\mapsto x(t)$ ($0\leq t\leq 1$) tel que $(x(t),e)$ soit une paire de Fredholm pour tout $t$. Enfin dans le paragraphe 6 on se restreint à la dimension finie pour montrer le lien entre cet indice et ceux de Clerc, Koufany et {\O}rsted.  

\vspace{0.7cm}

\noindent {\bf Remerciements}. Je tiens a remercier K.H. Neeb de m'avoir signaler une erreur dans la version pr\'ec\'edente de cet article.

\begin{nota}
Si $X$ est un espace topologique, On note $C(X)$ l'algèbre des fonctions complexes continues sur $X$. Si $E$ et $F$ sont deux espaces de Banach, on note $L(E,F)$ l'espace de Banach des opérateurs linéaires continus de $E$ dans $F$ muni de la norme d'opérateur et on pose $L(E)=L(E,E)$. Si $\mathcal{B}$ est une algèbre de Banach (associative) complexe et $x\in\mathcal{B}$, on note $sp(\mathcal{B},x)$ le spectre de $x$ dans $\mathcal{B}$. Lorsque $\mathcal{B}=L(E)$ on note simplement $sp(x)$ s'il n'y a pas d'ambiguïté.  
\end{nota}

\section{$JB^*$-triples et domaines bornés symétriques}\label{parnot}
 
Un $JB^*$-triple est la donnée d'un espace de Banach complexe $(E,\nm{.})$ et d'une application (on note $\overline{E}$ l'espace conjugué de $E$)
$$Q\ :\ E\rightarrow L{(\overline{E},E)}$$
quadratique et continue, telle que si l'on note
$$\{x,y,z\}=L(x,y)z=\frac12(Q(x+y)-Q(x)-Q(y))z$$
le système triple associé on ait l'{\it identité triple de Jordan} :
\begin{equation*}
\{u,v,\{x,y,z\}\}=\{\{u,v,x\},y,z\}-\{x,\{v,u,y\},z\}+\{x,y,\{u,v,z\}\}
\end{equation*}               
et les propriétés suivantes pour tout $x$ de $E$ :
\begin{enumerate}
\item 
$L(x,x)$ est un opérateur hermitien positif,
\item
$\nm{\{x,x,x\}}=\nm{x}^{3}$.
\end{enumerate}
Une algèbre de Jordan Banach est un espace de Banach $(E,\nm{.})$ muni d'un produit commutatif $x\circ y$ tel que
\begin{enumerate}
\item 
$\nm{x\circ y}\leq\nm{x}\nm{y},$
\item
$x\circ(x^2\circ y)=x^2\circ(x\circ y),\quad\forall x,y\in E$.
\end{enumerate}
Supposons $E$ complexe et muni d'une involution antilinéaire $*$. Alors
$$\{x,y,z\}=x\circ(y^*\circ z)+z\circ(y^*\circ x)-(x\circ z)\circ y^*$$
vérifie l'identité triple de Jordan et $E$ est appelée une $JB^*$-algèbre si l'on a 
$$\nm{\{x,x,x\}}=\nm{x}^{3}, \quad \forall x\in E.$$
Si $E$ possède un neutre, c'est alors un $JB^*$-triple.
Une algèbre de Jordan Banach réelle $A$ est appelée $JB$-algèbre si l'on a
\begin{enumerate}
\item
$\nm{x^2}=\nm{x}^2,$
\item
$\nm{x^2}\leq\nm{x^2+y^2},\quad \forall x,y\in A$.
\end{enumerate}
La partie réelle d'une $JB^*$-algèbre est une $JB$-algèbre et
réciproquement, étant donnée une $JB$-algèbre $A$, il existe sur $E=A\otimes\cc$ une unique norme prolongeant celle de $A$ et qui fait de $E$ (muni du produit étendu par linéarité) une
$JB^*$-algèbre (cf \cite{WR}).

Si $E$ est une $C^*$-algèbre de produit $xy$ alors $E$ muni du produit de Jordan 
$$x\circ y=\frac12(xy+yx)$$
devient une $JB^*$-algèbre.
Une $JB^*$-algèbre qui est isomorphe à une sous $JB^*$-algèbre (ie. un sous-espace fermé stable par le produit de Jordan) d'une $C^*$-algèbre est dite {\it spéciale}.

Un ouvert connexe et borné $\D$ d'un espace de Banach $E$ est appelé domaine borné symétrique si à chacun de ses points on peut associer un automorphisme holomorphe involutif de $\D$ dont il est un point fixe isolé. Un tel domaine est homogène sous son groupe d'automorphismes et biholomorphiquement équivalent à un domaine borné cerclé (ie. contenant l'origine et invariant sous l'action des nombres complexes de module 1) et étoilé par rapport à l'origine \cite{VI}. Une telle réalisation est unique à isomorphisme linéaire près (car un biholomorphisme d'un domaine cerclé conservant l'origine est linéaire). L'ensemble des champs de vecteurs complets sur $\D$ est une algèbre de Lie Banach et le groupe des biholomorphismes de $\D$ peut être muni d'une structure de groupe de Lie Banach réelle dont l'algèbre de Lie s'y identifie (cf. \cite{VI,UP1,UPsbmj}). Lorsque $\D$ est réalisé comme domaine cerclé cette algèbre de Lie que l'on notera $\g$ se décompose suivant les espaces propres de l'action de la symétrie à l'origine :
$$\g=\gk\oplus\gp$$
de sorte que $\gk$ est constitué de champs linéaires et que l'application
\begin{align*}
\gp&\rightarrow E\\ 
X&\mapsto X(0)
\end{align*}
est un isomorphisme de Banach. De plus il existe une application $Q\ :\ E\rightarrow L(\overline{E},E)$ quadratique et continue telle que pour tout $v\in E$ l'unique champ $X_v$ de $\gp$ tel que $X_v(0)=v$ s'écrive
$$X_v(z)=v-Q(z)v.$$
Cette application fait de $E$ un $JB^*$-triple dont la boule unité coïncide avec $\D$. Réciproquement la boule unité d'un $JB^*$-triple est un domaine borné symétrique (cf. \cite{KA1,KA2}).

Un élément $x$ d'un $JB^*$-triple $E$ est dit inversible si $Q(x)$ l'est. On note
$$x^\#=Q(x)^{-1}x.$$ 
On appelle tripotent tout élément tel que $Q(x)x=x$
et on note $\Sigma$ l'ensemble des tripotents inversibles de $E$. C'est une sous-variété banachique de $E$.  
Si $e\in\Sigma$ alors le produit
$$x\circ y:=L(x)y:=\{x,e,y\}$$ 
et l'involution $Q(e)$ font de $E$ une $JB^*$-algèbre de neutre $e$ que l'on notera $E^{(e)}$ et le système triple associé à $E^{(e)}$ est bien celui de $E$. Pour cette raison on appelle parfois $\Sigma$ l'ensemble des unités de $E$. On notera $A(e)$ la partie réelle de $E^{(e)}$ et 
$$P(x)=Q(x)Q(e)$$
la représentation quadratique. Un élément $x$ dans $E$ est donc inversible si et seulement si $P(x)$ l'est et on définit son inverse dans $E^{(e)}$ par
$$x^{-1}=P(x)^{-1}x=Q(e)x^{\#}.$$
Notons $x^*=Q(e)x$. Alors
$$\Sigma=\{x\in E\mid x^*=x^{-1}\}.$$

La notion d'inversibilit\'e dans une alg\`ebre de Jordan que nous avons introduite est due \`a N. Jacobson, qui a montr\'e qu'elle est équivalente à la définition classique : $x$ est inversible si et seulement si il existe un élément $y$ tel que $x\circ y=e$ et $x^2\circ y=x$, auquel cas $y$ est unique et est appelé l'inverse de $x$ (si E est une algèbre de Jordan  spéciale, alors ces deux propriétés sont équivalentes à $xy=yx=1$, cf. \cite[p.51]{JAsr}).

L'operateur de Bergman de $E$ est par d\'efini par 
$$B(x,y):=Id-2L(x,y)+Q(x)Q(y).$$
Le couple $(x,y)$ est dit \emph{transverse} lorsque $B(x,y)$ est inversible. Lorsque $y=e\in\Sigma$, on a $B(x,e)=Q(x-e)Q(e)=P(x-e)$. Donc le couple $(x,e)$ est transverse si et seulement si il est inversible. 

Le spectre de $x$ dans $E^{(e)}$, not\'e $Sp(x,e)$, est l'ensemble des nombres complexes $\lambda$ tels que $\lambda e-x$ n'est pas inversible. Alors d'apr\`es un th\'eor\`eme de J. Martinez Moreno (cf. \cite{MAR} et \cite{KA2}) : 
$$sp(L(x))\subset \frac12(Sp(x,e)+Sp(x,e)),$$ 
et
$$sp(P(x))\subset Sp(x,e)Sp(x,e).$$

On appelle tripotent régulier un tripotent $x$ tel que $\ker L(x,x)=0$ et on note $S$ leur ensemble. Lorsque la dimension de $E$ est finie (la théorie devient celle des systèmes triples de Jordan hermitiens positifs cf. \cite{LObsd}), si $\Sigma$ est non vide alors $\Sigma=S$ (car $S$ est homogène sous le groupe des automorphismes du système triple). En dimension infinie ce n'est plus le cas, mais $S$ s'identifie toujours à la frontière extrémale (au sens de la convexité) de $\overline{\D}$, et $\Sigma$ est une réunion de composantes connexes de $S$ (cf. \cite{KAUP,BKU}).

\section{La Lagrangienne comme frontière de Shilov de $Sym(H)$}

Soit $(H,\ps{\cdot,\cdot})$ un espace de Hilbert complexe 
(séparable). Le produit hilbertien $\ps{\cdot,\cdot}$ est 
antilinéaire par rapport à la seconde variable. Soit $\tau$ 
une involution (ie. une application $\cc$-antilinéaire 
involutive) isométrique de $H$. On note $\Sym(H)$
l'espace de Banach des opérateurs symétriques pour la forme 
bilinéaire symétrique 
$$(\cdot,\cdot)=\ps{\cdot,\tau(\cdot)}.$$
Pour $z\in\Li(H)$, on pose
$$\ov{z}=\tau\circ z\circ\tau.$$
L'espace de Banach $\Sym(H)$ muni du produit triple
$$\{x,y,z\}=\frac12(x\ov{y}z+z\ov{y}x)$$ est un $JB^*$
-triple. En effet, c'est un sous-système triple de Jordan 
fermé de $\Li(H)$ et on peut donc appliquer 
\cite[20.9]{UPsbmj}. 

Soit $z\in\Sym(H)$. On voit facilement que la notion d'inversibilit\'e co\"incide avec celle des op\'erateurs
(en effet, si $x$ est inversible comme op\'erateur alors $Q(x)$ l'est, et si $Q(x)$ est 
inversible, alors $\id=x(Q(x)^{-1}\id) x$
et donc $x$ est inversible),
et que, puisque,
$$\ov{x}x=\id \Rightarrow  x\ov{x}=\id,$$
les tripotents maximaux sont inversibles :
$$\Sigma=S=\{x\in \Sym(H)\mid \overline{x}x=\id\}.$$ 
L'opérateur de Bergman s'écrit
\begin{align*}
B(x,y)z&=z-(x\overline{y}z+z\overline{y}x)+x\overline{y}z\overline{y}x\\
&=(1-x\overline{y})z(1-\overline{y}x),
\end{align*}
et lorsque $x$ et $y$ sont dans $\Sigma$ on a
$$B(x,y)z=Q(y-x)Q(y)z=(1-xy^{-1})z(1-y^{-1}x),$$
et le couple $(x,y)$ est transverse si et seulement si $y-x$ est 
inversible. 

Considérons la structure de $JB^*$-algèbre sur $\Sym(H)$ 
définie par le tripotent inversible $\id$. Le produit s'écrit
$$x\circ y=\frac12(xy+yx)$$
et l'involution
$$x^*=\ov{x}.$$
Soit $H_0=\ker(\tau-\id)$ la forme réelle de $H$ associée à 
$\tau$. Alors la partie autoajointe de la 
$JB^*$-algèbre $\Sym(H)$ s'identifie à l'espace $\Sym(H_0)$ 
des opérateur symétriques de $H_0$ 
(qui est donc une $JB$-algèbre).

Introduisons maintenant un peu de vocabulaire et quelques 
notations. Soit $(\mathcal H,\ps{\cdot,\cdot})$ un espace de 
Hilbert 
réel ou complexe. Une forme bilinéaire antisymétrique 
$\omega$ continue  et fortement non-dégénérée (ie. telle que 
l'application $\mathcal H\ra \mathcal H'$, 
$\xi\mt\omega(\cdot,\xi)$ est bijective) est appelée 
\emph{forme symplectique}. Supposons $\mathcal H$ muni d'une 
telle forme. On dit alors que $\mathcal H$ est un \emph
{espace de Hilbert symplectique}.  
Pour un sous-espace $F\subset \mathcal H$, on note $F^\circ$ 
l'orthogonal de $F$ pour $\omega$, alors que l'on note 
$F^\perp$ l'orthogonal pour la structure Hilbertienne.
Un \emph{lagrangien} de $\mathcal H$ est un sous-espace  
$\lambda$ tel que $\lambda^\circ=\lambda$. Un lagrangien est 
automatiquement fermé (car $(F^\circ)^\circ=\ov{F}$ pour 
tout sous espace $F$).
On appelle \emph{Lagrangienne} l'ensemble des lagrangiens de 
$\mathcal H$, et on la note 
$\Lambda(\mathcal H)$.

On pose $\mathbb{H}=H\oplus H=\{\eta\oplus \xi \mid 
\xi,\eta\in H\}$. C'est un espace de Hilbert pour 
la forme hermitienne
$$\left\langle \eta\oplus \xi,\eta'\oplus \xi'\right\rangle 
=\left\langle \eta,\eta' 
\right\rangle+\left\langle \xi,\xi'\right\rangle,$$
et on  muni $\hh$ d'une structure symplectique (complexe) en 
posant 
$$\omega(\eta\oplus \xi,\eta'\oplus 
\xi')=(\eta,\xi')-(\xi,\eta').$$
L'involution $\tau$ s'\'etend \`a $\hh$ en posant 
$$\tau(\eta\oplus \xi)=\tau(\eta)\oplus\tau(\xi).$$
Alors $\hh_0=H_0\oplus H_0$ est la forme réelle de 
$\hh$ asociée à $\tau$, et puisque la  forme symplectique v\'erifie 
$$\omega(\tau(\eta\oplus \xi),\tau(\eta'\oplus \xi'))=\overline{\omega(\eta\oplus 
\xi,\eta'\oplus \xi')},$$
on peut la restreindre \`a $\hh_0$ en une forme symplectique r\'eelle. Nous allons 
montrer comment l'ensemble $\Sigma$ s'identifie à la lagrangienne $\Lambda(\hh_0)$ de $\hh_0$. 

Commençons par envoyer $\Sym(H)$ dans $\Lambda(\hh)$. 

Notons $H_1=\ho$ et $H_2=\oh$, $\pi_1$ et $\pi_2$ les projections sur $H_1$ (resp. $H_2$) 
parall\`element \`a $H_2$ (resp. $H_1$).
Si $x\in \Li(H)$ alors 
$$G(x):=\{x\xi\oplus \xi\mid \xi\in H\}$$
est un sous-espace ferm\'e de 
$\hh$. Il est de plus transverse \`a
$H_1$ (ie. $G(x)\oplus H_1=\hh$) car $\xi\oplus \eta=\xi'\oplus(x\xi'+\eta')$, 
$\xi,\eta,\xi',\eta' \in H$ se r\'esout de mani\`ere unique en $\xi=\xi'$, $\eta'=\eta-x\xi$. 
R\'eciproquement, soit $F$ un sous-espace ferm\'e et transverse \`a $H_1$ et 
$\pi:F\rightarrow H_2$ la restriction de $\pi_2$ \`a $F$. L'application $\pi$ est bijective 
parce que $F$ est transverse et comme $H_1$ est un suppl\'ementaire ferm\'e $\pi$ est 
continue et donc 
d'apr\`es le th\'eor\`eme de Banach elle est bicontinue. Alors $\pi_1\circ\pi^{-1}\in \Li(H)$ 
et $G(\pi_1\circ\pi^{-1})=F$.
Remarquons que $H_1$ et $H_2$ 
sont dans $\Lambda(\hh)$.  
Si $x\in \Li(H)$ on note $\trans{x}$ le transpos\'e de $x$ par rapport \`a $(.,.)$. Alors 
$G(\trans{x})=G(x)^\circ$. En effet l'inclusion $G(\trans{x})\subset G(x)^\circ$ est clair et 
si il n'y avait pas \'egalit\'e on
aurait $G(x)^\circ\cap H_2\neq\{0\}$ ce qui impliquerait $H_1\cap      H_2\neq\{0\}$. 
L'application $G$ induit donc une bijection entre $\Sym(H)$ et les 
lagrangiens transverses \`a $H_1$. 

Posons $J(\eta\oplus \xi)=(-\xi)\oplus \eta$. Alors 
$\omega(\cdot,\cdot)=\left(J\cdot,\cdot\right)$ et pour tout $\lambda\in\Lambda(\hh)$, 
$\lambda^\perp=J\tau(\lambda)$.

Pour caract\'eriser l'image de $\Sigma$ par l'application $G$ introduisons la forme 
hermitienne 
$$h(\eta\oplus \xi,\eta'\oplus \xi')=\langle \xi,\xi'\rangle-\langle \eta,\eta'\rangle.$$
\begin{pro}
Soit $\lambda$ un lagrangien sur lequel $h$ est une forme positive. Alors $\lambda$
est transverse \`a $H_1$.
\end{pro}
\begin{proof}
Soit $\lambda$ un lagrangien sur lequel $h$ est une forme positive. Si $\eta\oplus 0\in\lambda$ alors 
$$h(\eta\oplus 0,\eta\oplus 0)=-\langle \eta,\eta\rangle\geq0$$
donc $\eta=0$ et $\lambda\cap H_1=\{0\}$. 
Comme $$(\lambda+H_1)^{\perp}=\lambda^{\perp}\cap{H_1}^{\perp}=J\tau(\lambda)\cap 
J\tau(H_1)=J\tau(\lambda\cap{H_1})=\{0\},$$ 
il suffit de montrer que $\lambda+H_1$ est ferm\'e. Soit $(\zeta_{n})_{\nn}$ une suite de 
$\lambda+H_1$ qui converge vers $\zeta\in\hh$. Pour tout entier $n$, 
$\zeta_{n}=\xi_{n}+\eta_{n}+\eta'_{n}$ avec 
$\xi_{n}\in H_2,\ \eta_{n},\eta'_{n}\in H_1\ \textrm{et}\ \xi_{n}+\eta_{n}\in\lambda$.
$\xi_{n}$ est la projection orthogonale de $\zeta_{n}$ sur $H_2$ et converge
donc vers la projection orthogonale $\xi$ de $\eta$ sur $H_2$. D'autre part, 
$h(\xi_{n}+\eta_{n},\xi_{n}+\eta_{n})=\langle \xi_{n},\xi_{n}\rangle-\langle 
\eta_{n},\eta_{n}\rangle\geq0$ donc $(\eta_{n})_{\nn}$ est born\'ee et on peut extraire 
une suite, toujours not\'ee $(\eta_{n})_{\nn}$, qui converge faiblement vers $\eta$. Mais 
$H_1$ est ferm\'e pour la topologie forte et convexe donc ferm\'e 
pour la topologie faible et donc $\eta\in H_1$. Comme $\xi_{n}+\eta_{n}$ converge
faiblement vers $\xi+\eta$ et que $\eta'_{n}$ converge faiblement vers $\eta'=\zeta-\xi-\eta$
, on en d\'eduit de m\^eme que $\xi+\eta\in\lambda$ et $\eta'\in H_1
$. Par unicit\'e de la limite on obtient la d\'ecomposition $\zeta=\xi+\eta+\eta'$ qui nous 
permet de conclure que $\zeta\in\lambda+H_1$. Finalement, $\lambda$ est bien transverse \`a 
$H_1$.
\end{proof}
\noindent Alors $h$ s'annule sur $\lambda=G(x)$ si et seulement si pour tout $\xi\in H$, 
 $\left\langle x\xi,x\xi\right\rangle=\left\langle\xi,\xi\right\rangle$, ie. si et seulement si 
$\left\langle(1-x^*x)\xi,\xi\right\rangle=0$ ce qui \'equivaut par polarisation \`a $x^*x=1$. En résumé,
\begin{align*}
\Sym(H)&\overset{G}{\hookrightarrow}\Lambda(\hh)\\
\Sigma&\simeq\{\lambda\in\Lambda(\hh)\mid h_{\mid \lambda\times\lambda}=0\}.
\end{align*}

On d\'efinit la transform\'ee de Cayley sur $\hh$ par
$$C(\eta\oplus \xi)=\frac{1}{\sqrt{2}}((\eta+i\xi)\oplus(i\eta+\xi)).$$
On a 
$$\omega(C\cdot,C\cdot)=\omega(\cdot,\cdot)\quad\text{et}\quad
 ih(\cdot,\cdot)=\omega(C\cdot,\tau(C\cdot)).$$
Donc $C$ conserve les lagrangiens de $\hh$, et les 
lagrangiens sur lesquels $h$ s'annule sont 
transform\'es en les lagrangiens stables par $\tau$. Il ne 
reste plus qu'à identifier l'ensemble des lagrangiens 
stables par $\tau$ et $\Lambda(\hh_0)$. 

Si $F_0$ est un sous-espace de $\hh_0$, alors $F_0\oplus 
iF_0$ est un sous-espace complexe de $\hh$ stable par $\tau$
. Réciproquement, si $F$ est un sous-espace de $\hh$ stable 
par $\tau$, alors $F=F\cap\hh_0\oplus F\cap 
i\hh_0=F\cap\hh_0\oplus i(F\cap\hh_0)$. De plus il est clair 
que si $F_0$ est un lagrangien réel, alors $F_0\oplus 
iF_0$ est un lagrangien complexe et que si $F$ est un 
lagrangien complexe, alors $F\cap\hh_0$ est un lagrangien réel. On a donc une bijection entre la Lagrangienne r\'eelle et l'ensemble $\Lambda(\hh)^{\tau}$ des 
lagrangiens complexes stables par $\tau$.
Finalement on a bien une bijection entre $\Sigma$
et la Lagrangienne r\'eelle:
$$\Sigma\overset{G}{\simeq}\{\lambda\in\Lambda(\hh)\mid h_{\mid \lambda\times\lambda}=0
\}\overset{C}{\simeq}\Lambda(\hh)^{\tau}\simeq\Lambda(\hh_0).$$

La Lagrangienne $\Lambda(\hh_0)$ est munie d'une structure de 
variété banachique (cf. \cite{FU}). La bijection que nous 
avons décrite est alors un difféomorphisme. Nous n'écrivons 
pas les détails.  

Réciproquement, partons maintenant d'un espace de Hilbert 
symplectique réel 
$(\hh_0,\omega,\langle\cdot,\cdot\rangle)$.
On peut supposer, quitte à remplacer le produit scalaire par 
un autre définissant une norme équivalente, que la forme 
symplectique et le produit scalaire sont \emph{compatibles}, 
c'est-\`a-dire qu'ils sont li\'es par la relation 
$$\omega(\cdot,\cdot)=\langle 
J\cdot,\cdot\rangle$$
o\`u $J$ est \`a la fois un op\'erateur 
orthogonal et une structure 
complexe (cf. \cite[Appendix D]{FU}). 
Soit $H_0\in\Lambda(\hh_0)$. Alors $H_0$ 
est ferm\'e et $H_0^{\perp}=JH_0$. 
Suivant la d\'ecomposition 
$\hh_0=H_0\oplus JH_0\simeq H_0\oplus H_0$,
$$\omega(\eta\oplus \xi,\eta'\oplus \xi')=\langle 
\eta,\xi'\rangle-\langle \xi,\eta'\rangle$$
et on peut \'etendre $\omega$ et 
$\left\langle\cdot,\cdot\right\rangle$ \`a 
$\hh=\hh_0\otimes\mathbb{C}$ et poser $H=H_0\oplus iH_0$ 
pour se retrouver dans la situation du paragraphe précédent.

La notion de \emph{paire de Fredholm} est essentielle dans 
la définition de l'indice de Maslov en dimension infinie.
La paire de lagrangiens $(\lambda,\mu)\in \Lambda(\hh_0)^2$ 
est appel\'ee paire de Fredholm si
$$\dim \lambda\cap\mu<\infty\quad\text{et}\quad\dim 
\hh_0/(\lambda+\mu)<\infty.$$ 
La Fredholm-Lagrangienne relativement à $\lambda$ est alors
$$\F\Lambda_\lambda=\{\mu\in\Lambda(\hh_0)\mid (\mu,\lambda)\ 
\text{est une paire de Fredholm}\}.$$
L'indice de Maslov relativement à $\lambda$ est défini pour 
les chemins (continus) dans la Fredholm-Lagrangienne 
relativement à $\lambda$ (cf. \cite{BOFU,FU}). Nous voulons 
maintenant traduire cette notion dans la réalisation de la 
Lagrangienne
comme ensemble des tripotents inversibles de $\Sym(H)$. 
Soient $x$ et $y$ deux op\'erateurs de $H$, et $G(x)$ et 
$G(y)$ leurs graphes dans $\hh$. Alors
$$\ker(y-x)=\pi_2(G(x)\cap G(y)),$$
et
\begin{equation*}
H/\ker(y-x)\simeq G(x)/G(x)\cap G(y)
\end{equation*}

\noindent D'autre part, 
\begin{align*}
G(x)+G(y)&=\{x\xi\oplus \xi+y\xi'\oplus \xi'\mid \xi,\xi'\in H\}\\
&=\{(x\xi+y\xi')\oplus(\xi+\xi')\mid\xi,\xi'\in H\}\\
&=\{(x\zeta+(y-x)\xi')\oplus\zeta\mid\zeta,\xi'\in H\}\\
&=G(x)+((y-x)H\oplus0).
\end{align*}
En consid\'erant l'application
\begin{gather*}
\hh\rightarrow H_1\rightarrow H_1/((y-x)H\oplus0)\\
\xi\oplus\eta\mapsto 
(\eta-x\xi)\oplus0\mapsto((\eta-x\xi)\oplus0)+((y-x)H\oplus0) 
\end{gather*}
on obtient l'isomorphisme
\begin{equation*}
\hh/G(x)+G(y)\simeq H/(y-x)H.
\end{equation*}
Supposons maintenant que $x$ et $y$ sont dans $\Sigma$, et 
soient $\lambda$ et $\mu$ les lagrangiens associées.
Comme $C$ et l'application qui à un lagrangiens réel associe 
le lagrangien complexe qu'il engendre, respectent 
l'intersection et la somme, on en 
d\'eduit que $\lambda$ et $\mu$ forment une paire transverse (ie. $\lambda\oplus\mu=\hh_0$) si et seulement si 
$y-x$ est inversible, ie. si et seulement si $(x,y)$ est transverse, et forment une paire de 
Fredholm si et seulement si $y-x$ est un op\'erateur de Fredholm sur $H$. De plus, $\pi_2$ 
étant injective sur $G(x)\cap G(y)$, on a 
$$\dim_{\cc}\ker(y-x)=\dim\lambda\cap\mu.$$ 

On peut considérer $\hh_0$ comme un espace de Hilbert 
complexe grâce à la la 
structure presque complexe $J$ et au produit hilbertien
$$\ps{\cdot,\cdot}_J=\ps{\cdot,\cdot}-i\omega(\cdot,\cdot).$$
On note alors $U(\hh_0,J)$ le 
groupe des opérateur unitaires. Puisque 
$\omega(\cdot,\cdot)=\ps{J\cdot,\cdot}$, on voit
que $U(\hh_0,J)$ agit sur $\Lambda(\hh_0)$. Cette action est 
de plus transitive. En effet, soient $\lambda$ et $\mu$ deux 
lagrangiens. Soit $(e_n)_{n\in\nn}$ une base hilbertienne 
réelle de $\hh_0$. Alors comme 
$\hh_0=\lambda\oplus\lambda^{\perp}=\lambda\oplus J\lambda$, 
$(e_n)_{n\in\nn}$ est une base hilbertienne complexe de 
$\hh_0$. De même  
 une base hilbertienne réelle $(e'_n)_{n\in\nn}$ de 
$\mu$ est une base hilbertienne complexe de 
$\hh_0$. On sait qu'il existe un opérateur unitaire envoyant 
la base $(e_n)_{n\in\nn}$ sur la base $(e'_n)_{n\in\nn}$, 
et un tel opérateur envoie $\lambda$ sur $\mu$.

Nous voulons maintenant transporter l'action de $U(\hh_0,J)$ 
en une action d'un groupe (à caractériser) sur $\Sigma$.
Soit $U\in U(\hh_0,J)$. Notons $U_\cc$ l'extension 
$\cc$-linéaire de $U$ à $\hh$. Alors $U_\cc$ agit sur $\Lambda(\hh)$, et préserve $\Lambda(\hh)^\tau$. Pour tout $a,b,c,d\in H_0$, un calcul montre que
$$C^{-1}U_\cc C((a+ib)\oplus(c+id))=(a'-ib')\oplus(c''+id''),$$
où $a',b',c'',d''\in H_0$ sont définis par
$$a'\oplus b'=U(a\oplus(-b))\quad \text{et}\quad c''\oplus 
d''=U(c\oplus d).$$
Remarquons que ${C^{-1}U_\cc C}$ laisse $H_1=H\oplus0$ et 
$H_2=0\oplus H$ stables et notons 
$$u={(C^{-1}U_\cc C)}_{\mid H_1}\quad\text{et}\quad 
v={(C^{-1}U_\cc C)}_{\mid H_2}.$$
On considère $u$ et $v$ comme des opérateurs sur $H$. Ce 
sont des opérateurs unitaires de $H$ car $C^{-1}U_\cc C$ est 
unitaire.
Soit $x\in \Sym(H)$. Alors 
$$C^{-1}U_\cc C(G(x))=G(uxv^{-1}).$$
Montrons que $v^{-1}=\trans{u}$. Il suffit de montrer que pour tout $a,b\in H_0$, 
$$(u(a+ib),v(a+ib))=(a+ib,a+ib),$$
donc que
$$(a'-ib',a''+ib'')=(a+ib,a+ib),$$
ou encore, en développant (le produit hermitiens et la forme 
bilinéaire coïncidant sur $H_0$), que
$$\ps{a', a''}+\ps{b', b''}+i(\ps{a', b''}-\ps{b', a''})=
\ps{a,a}-\ps{b,b}+2i\ps{a,b}.$$
Mais comme $U\in U(\hh_0,J)$,
$$\ps{U(a\oplus(-b)),U(a\oplus b)}_J=\ps{a\oplus(-b),a\oplus 
b}_J,$$
donc
$$\ps{a'\oplus b',a''\oplus b''}_J=\ps{a\oplus(-b),a\oplus 
b}_J,$$
ce qui donne la relation voulue. 
En résumé, l'action de $U(\hh_0,J)$ sur $\Lambda(\hh_0)$ se 
transporte en une action (transitive) du groupe unitaire 
$U(H)$ sur $\Sigma$ (on peut en 
effet montrer que l'on obtient bien tout $U(H)$), et cette 
action est la restriction de l'action de $U(H)$ sur $\Sym(H)$
définie par 
$$U(H)\times\Sym(H)\ra\Sym(H),\quad (u,z)\mt u z \trans{u}.$$
Remarquons pour conclure que $U(H)$ agit par automorphismes du système triple de Jordan $\Sym(H)$.

\section{Les paires de Fredholm et l'indice de transversalité}

Dans cette partie on considère un $JB^*$-triple $E$ tel que $\Sigma$ n'est pas vide. 
Soit $e\in\Sigma$. Les notations sont celles du paragraphes 2.
Soit $x\in\Sigma$ et soit $C^*(x,e)$ la sous-algèbre 
fermée de $E^{(e)}$ engendrée 
par $e$, $x$ et $x^*=Q(e)x$.
\begin{pro}
Soient $(x,e)\in\Sigma^2$. Alors
\begin{enumerate}[(i)]
\item $C^*(x,e)$ est associative et c'est donc une $C^{*}$-alg\`ebre commutative.
\item Le spectre $U_{x,e}$ de $x$ dans $C^*(x,e)$ est contenu dans le cercle unité, et c'est 
aussi le spectre de $x$ dans $E^{(e)}$.
\item La paire $(x,e)$ est transverse, ie. $B(x,e)$ est inversible, si et seulement si $1\not\in U_{x,e}$.
\end{enumerate}
\end{pro}
\begin{proof}
Dans $E^{(e)}$, $x^*=x^{-1}$. Or
on a $[L(x),L(x^{-1})]=0$ (cf. \cite[19.26]{UPsbmj}) et donc 
$C^*(x,e)$ est fortement associative, en particulier 
associative. Le système triple $C^*(x,e)$ est donc lui aussi 
associatif et l'on a (cf. \cite[20.32]{UPsbmj}), pour tout 
$u,v\in C^*(x,e)$, $\nm{u\circ v}\leq\nm{u}\nm{v}$. On écrit alors 
comme dans \cite[20.33]{UPsbmj}, pour $z\in C^*(x,e)$,
$$\nm{z}^3=\nm{\{z,z,z\}}=\nm{z\circ(z^*\circ 
z)}\leq\nm{z}\nm{z^*\circ z}\leq\nm{z}^2\nm{z^*}=\nm{z}^3.$$
Donc $C^*(x,e)$ est une $C^*$-algèbre. Comme $x$ est 
unitaire dans cette $C^*$-algèbre, son spectre est contenu 
dans le cercle unit\'e. Or $C^*(x,e)$ est contenue dans 
une sous-algèbre fortement 
associative maximale (et fermée) de $E^{(e)}$, et le spectre de $x$ 
dans cette sous-algèbre est égal au 
spectre de $x$ dans $C^*(x,e)$, puisque celui-ci est égal à 
sa 
frontière. Mais cette 
sous-algèbre fortement associative maximale est pleine dans 
l'algèbre de Jordan $E^{(e)}$ (ie. ses \'el\'ements y sont inversible si et seulement si ils sont inversibles dans $E^{(e)}$, cf. \cite{HES, MAth}), et donc le spectre de $x$ 
dans $C^*(x,e)$ est égal au spectre de $x$ dans $E^{(e)}$. 
La dernière assertion découle immédiatement de la précédente 
et du fait que $B(x,e)=Q(x-e)Q(e)=P(x-e)$.
\end{proof}

Puisque $C^*(x,e)$ est engendrée (comme $C^*$-algèbre) par
$x$ et $e$,
on a l'isomorphisme de Gelf'and :
\begin{align*}
\mathcal{G}_{x,e}:&C^*(x,e)\rightarrow C(U_{x,e}),\\
&y\mapsto\hat{y},
\end{align*}
qui \`a l'élément $x$ associe la fonction $\hat{x}(\mu)=\mu$.

\begin{pro} Soit $(x,e)\in\Sigma^2$. On a
$$U_{e,x}=\ov{U_{x,e}}.$$
\end{pro}
\demo
Cela résulte du fait que $Q(e-\lambda x)$ est inversible si et seulement si $Q(x-\lambda^{-1} e)$ l'est.
\fdemo
Supposons maintenant que $1\not\in U_{x,e}$ ou bien que $1$ 
est isolé dans $U_{x,e}$. Alors la fonction caractéristique 
$\chi_{\{1\}}$ de $\{1\}$
est continue sur $U_{x,e}$. On note alors 
$$p=p(x,e)={\mathcal{G}_{x,e}}^{-1}(\chi_{\{1\}})$$
le projecteur associé à $1$, et
$$A_p(e)=\{p,A(e),p\}.$$ 
On dit que $1$ est de multiplicité finie si $A_p(e)$ est une 
$JB$-algèbre de rang finie, ie.
$$A_p(e)=A_1\oplus\dots\oplus A_q$$ 
où chaque $A_j$ est une algèbre de Jordan euclidienne simple 
(cf. \cite{FAKO}) ou un facteur spin (ie. la $JB$-algèbre, de 
rang $2$, $H\oplus\rr$ où $H$ est un espace de Hilbert),
le rang de $A_p(e)$ étant alors par définition 
$$\rk{A_p(e)}=\rk{(A_1)}+\dots+\rk{(A_q)}.$$
\begin{defi}
Soit $(x,e)\in\Sigma$. On dit que $(x,e)$ est une paire de 
Fredholm lorsque $(x,e)$ est transverse, ou lorsque $1$ est 
isolé dans $U_{x,e}$, et est de multiplicité finie. 
\end{defi}
\noindent On 
définit alors l'indice
de transversalité de la paire de Fredholm comme le rang de $A_p(e)$,
$$\mu(x,e)=\rk{A_p(e)}.$$
Lorsque $1$ est isolé mais que $A_p(e)$ n'est pas de rang 
fini, on pose $\mu(x,e)=\infty$.

\begin{ex}
Considérons le $JB^*$-triple $E=\Sym(H)$. Les notations sont 
celles du paragraphe précédent. En particuliers $\tau$ désigne 
l'involution de $H$ et $H_0$ la forme réelle associée.
Soit $(x,e)\in\Sigma^2$. Alors 
$xe^{-1}$ est unitaire, donc normal. Soit 
$C^*(xe^{-1})$ la sous-algèbre fermée de $\Li(H)$ 
engendrée par $\id$, $xe^{-1}$, et $(xe^{-1})^*=ex^{-1}$. 
Alors la multiplication à droite par $e$ est un isomorphisme 
de $C^*(xe^{-1})$ sur $C^*(x,e)$ (qui envoie $\id$ sur 
$e$ et $xe^{-1}$ sur $x$). Supposons que $1$ est isolé dans 
$U_{x,e}$. Notons $p$ le projecteur  associé à $1$ dans 
$C^*(x,e)$ et $p'=pe^{-1}$ le projecteur associé à $1$ dans 
$C^*(xe^{-1})$. L'opérateur $p'$ est une projection au sens 
usuel, et l'on a 
$$\ker(\id-xe^{-1})=p'H.$$
Nous avons vu au chapitre précédent que l'action du groupe 
unitaire $U(H)$ sur $\Sigma$ par automorphismes du système 
triple $E$ (définie par $z\mt vz\trans{v}$) 
est transitive. En particuliers, il existe $u\in U(H)$ tel 
que $e=u\trans{u}$. Alors $A(e)$ (la partie autoadjointe de 
$E$ pour l'involution définie par $e$) est isomorphe à 
$\Sym(H_0)$ 
: 
$$A(e)=A(u\id\trans{u})=uA(\id)\trans{u}\simeq 
A(\id)\simeq\Sym(H_0).$$
Soit $p''=u^{-1}p\,\trans{u}^{-1}$. C'est un projecteur de 
$\Sym(H_0)$, ie. une projection de $H$ laissant $H_0$ 
stable. De 
plus
$$A_{p''}(\id)=p''\Sym(H_0)p''\simeq\Sym(p''H_0),$$
comme on peut le voir en écrivant la "matrice" d'un opérateur $z\in \Sym(H_0)$ relativement à la décomposition $H_0=p''H_0\oplus(1-p'')H_0$,
et donc
$$A_{p}(e)\simeq\Sym(p''H_0).$$
Ainsi $A_{p}(e)$ est de rang fini si et seulement si $p''H_0$
est de dimension finie. De plus
$$\ker(\id-xe^{-1})=pe^{-1}H=pH=p\,\trans{u}^{-1}H=up''H,$$
donc 
$$\dim\ker(\id-xe^{-1})=\dim_\cc p''H=\dim_\rr p''H_0,$$
et lorsque l'un des deux membres est fini,
$$\dim\ker(e-x)=\rang{A_p(e)}.$$
Montrons que dans ce cas, $x-e$ est un opérateur de 
Fredholm. Puisque 
$$(\id-xe^{-1})^*=\id-ex^{-1}=(x-e)x^{-1}=(xe^{-1}-\id)ex^{-1}
,$$
on a
$$\cd\ov{\im(\id-xe^{-1})}=\dim\ker(\id-xe^{-1}),$$
et il reste à montrer que $\id-xe^{-1}$ est d'image fermée. 
Grâce à l'isomorphisme de Gelf'and on voit que 
$(1-(xe^{-1}-p'))$ est inversible, et l'on a
$$\id-p'=(\id-xe^{-1})(1-(xe^{-1}-p'))^{-1}.$$
Donc
$$\im(\id-xe^{-1})=\im(\id-p')=\ker p'$$   
est bien fermé.
Réciproquement, si $x-e$ est un opérateur de Fredholm, 
c'est aussi le cas de $\id-xe^{-1}$. 
Mais lorsque $0$ est isol\'e dans la 
frontière du spectre d'un opérateur de Fredholm, il est en fait isolé. 
Il en r\'esulte que $1$ est isolé dans $U_{x,e}$.
\end{ex}

Revenons au cas général. 

\begin{pro}
Soit $(x,e)\in\Sigma^2$. Alors $1$ est isolé dans $U_{x,e}$ 
si et seulement si $0$ est isolé dans $B(x,e)$.
\end{pro}
\begin{proof}
Posons $U=U_{x,e}=Sp(x,e)$. Dans l'algèbre de Jordan unitaire $E^{(e)}$ on a $B(x,e)=P(x-e)$, et le théorème de J. Martinez Moreno, 
$$sp(P(x-e))\subset (1-U)(1-U).$$
Donc si $1$ est isolé dans $U$, alors $0$ est isolé dans $sp(B(x,e))$. 

Pour établir la réciproque, on montre que 
$$\{{(1-\lambda)}^2 \mid \lambda\in U\}\subset sp(P(x-e)).$$

Commençons par monter que 
$$\{{(1-\lambda)}^2 \mid \lambda\in U\}\subset\partial_{ext}(1-U)(1-U),$$
où $\partial_{ext}K$ est la frontière de la composante connexe non bornée du complémentaire du compact K. Tout élément  $1-\lambda\in 1-\uu$ s'écrit de manière unique 
$2\cos{\frac\Theta2}e^{i\frac\Theta2}$ avec 
$-\pi<\Theta\leq\pi$, et alors 
$$(1-\lambda)^2=4\cos^2{\frac\Theta2}e^{i\Theta}.$$
Soient $2\cos{\frac\theta2}e^{i\frac\theta2}$ et 
$2\cos(\Theta-\frac{\theta}2)e^{i(\Theta-\frac{\theta}2)}$ 
dans $1-\uu$ : leur produit vaut
$$4\cos{\frac\theta2}\cos(\Theta-\frac{\theta}2)e^{i\Theta}=
2\big(
\cos{\Theta}+\cos(\Theta-\theta)\big)e^{i\Theta}.$$
Or la fonction $\theta\mt\cos{\Theta}+\cos(\Theta-\theta)$ 
est maximale pour $\Theta=\theta$. La demi-droite 
$]4\cos^2{\frac\Theta2},+\infty[e^{i\Theta}$ est donc 
entièrement contenue dans le complémentaire de $(1-U)(1-U)$, et cela 
implique notre assertion.
Considérons 
$$\mathcal{B}=\{T\in \Li(E)\mid TC^*(x,e)\subset C^*(x,e)\}.$$ 
C'est une sous-algèbre fermée de $\Li(E)$ qui contient $P(x-e)$.
Le spectre de $P(x-e)$ dans $\mathcal B$ est constitué de 
$sp(P(x-e))$ et, éventuellement, de certains de ses 
trous (ie. les composantes connexes 
bornées de son complémentaire). De plus, en considérant le morphisme 
$$\mathcal B\ra\Li(C^*(x,e)),\ T\mt T_{|C^*(x,e)},$$ 
on a, puisque $sp(P(x-e)_{|C^*(x,e)})=\{{(1-\lambda)}^2 \mid \lambda\in 
U\}$, l'inclusion
$$\{{(1-\lambda)}^2 \mid \lambda\in 
U\}\subset sp(P(x-e),\mathcal{B}).$$
Il résulte alors de
$$\{{(1-\lambda)}^2 \mid \lambda\in U\}\subset\partial_{ext}(1-U)(1-U)$$
que 
$$\{{(1-\lambda)}^2 \mid \lambda\in U\}\subset\partial sp(P(x-e),\mathcal{B}).$$ 
Comme
$\partial sp(P(x-e),\mathcal{B})\subset\partial sp(P(x-e))$, il vient 
$$\{{(1-\lambda)}^2 \mid \lambda\in U\}\subset\partial sp(P(x-e)).$$
Et donc 
si $1\in U$ n'est pas isolé, alors $0\in sp(P(x-e))$ n'est pas isolé.
\end{proof}





On note $\Sigma_e$ le composante connexe de $\Sigma$ contenant $e$, et $\F\Sigma_e$ 
l'ensemble des $x\in \Sigma$ tels que $(x,e)$ est une paire de Fredholm.
\begin{pro}
Soit $\Sigma_e$ la composante connexe de $\Sigma$ contenant $e$. Alors
$$\mathcal{F}\Sigma_e\subset\Sigma_e.$$
\end{pro}
\begin{proof}
Soit $x\in\mathcal{F}\Sigma_e$. Alors $U_{x,e}$ n'est pas $\mathbb{U}$ tout entier, puisque soit $1$ n'y est pas, soit $1$ y est mais est isolé, et on 
peut donc d\'efinir $\log x\in C^*(x,e)$, $\log$ \'etant une d\'etermination ad\'equate du 
logarithme. Alors $P(\exp{(\frac12\log x))e=x}$ et donc on a bien $x\in\Sigma_e$.
\end{proof}

Considérons $(x,e)\in\Sigma^2$ et soit $e^{i\theta}\in \uu$. On a 
$U_{x,e^{i\theta}e}=e^{-i\theta}U_{x,e}$, et donc si 
$e^{i\theta}$ est isolé dans $U_{x,e}$, alors $1$ est isolé dans $U_{x,e^{i\theta}e}$, et on peut définir $p(x,e^{i\theta}e)$ et $\mu(x,e^{i\theta}e)$. 

Un sous-ensemble $\sigma$ de $U_{x,e}$ est dit spectral lorsque il est à 
la fois ouvert et fermé. Cela revient à dire que la fonction caractéristique $\chi_\sigma$ 
est continue.

\begin{lem}
Soient $(x,e)\in\Sigma^2$ et $\sigma$ un sous-ensemble spectral de $U_{x,e}$. Si 
$p=\mathcal{G}_{(x,e)}^{-1}(\chi_\sigma)$ alors $P(p)x$ et $p$ sont deux unit\'es de $P(p)E$ et
on a les propri\'et\'es suivantes :
\begin{enumerate}[(i)]
\item
$P(p)C^*(x,e)=C^*(P(p)x,p)\subset C^*(x,e)$,
\item
$U_{P(p)x,p}=\sigma$.
\end{enumerate}
\end{lem}
\begin{proof}
Comme $P(p)E$ est un sous syst\`eme triple de $E$, il est clair que $p$ et $P(p)$ sont des tripotents de $P(p)E$. Pour voir que $P(p)x$ y est inversible, on montre que $P(p)x^{-1}$ est l'inverse de $P(p)x$ dans l'alg\`ebre de Jordan $P(p)E^{(e)}$. Mais en calculant dans $C(U_{x,e})$ on voit facilement que $\{P(p)x,p,P(p)x^{-1}\}=p$ et que $\{(P(p)x)^2,p,P(p)x^{-1}\}=P(p)x$.\\
\noindent(1) Puisque $p\in C^*(x,e)$, on a
$$P(p)x^2=pxxp=pxxp^2=pxpxp=pxppxp=(P(p)x)^2.$$
(2) Si $\lambda\not\in U_{P(p)x,p}$ alors il existe $z\in C^*(P(p)x,p)$ tel que 
$$(\lambda p-P(p)x)z=p.$$ 
Soit $g(\mu)=(1-\chi_\sigma(\mu))(\lambda-\mu)^{-1}$ pour $\mu\in U_{P(p)x,p}$. Alors 
$$g(x)(\lambda e-x)=e-p.$$ 
Gr\^ace \`a (1), $(\lambda e-x)z=(\lambda e-x)P(p)z=(\lambda e-x)pzp=(\lambda p-P(p)x)z=p$ et 
$$(g(x)+z)(\lambda e-x)=e-p+p=e.$$
Donc $\lambda\not\in U_{x,e}$. R\'eciproquement, si $\lambda\not\in\sigma$ soit $h(\mu)=\chi_\sigma(\mu)(\lambda-\mu)^{-1}$ pour $\mu\in U_{x,e}$. Alors $h(x)(\lambda e-x)=p$ donc 
$$P(p)h(x)(\lambda p-P(p)x)=ph(x)p((\lambda p-pxp)=ph(x)(\lambda e-x)p=p.$$
\end{proof}

\noindent Pour $0<\varepsilon<\pi$ on note 
$$\mathcal 
A_\varepsilon=\{e^{i\theta}\mid 0<\nm{\theta}\leq \varepsilon\}.$$

\begin{lem}[Perturbation de l'indice de transversalité]
\label{pert}
Soit $(x,e)$ une paire de Fredholm.
Il existe $0<\varepsilon<\pi$ tel que $1$ est la seule 
valeur spectrale de $x$ dans $\mathcal{A}_\varepsilon$. 
Il existe un voisinage $\mathcal V$ de $x$ tel que pour tout 
tripotent inversible $y\in\mathcal V$, le spectre de $y$ dans 
$\mathcal{A}_\varepsilon$ est fini et ne contient pas $e^{\pm i\varepsilon}$, et
$$\mu(x,e)=\sum_{\nm{\theta}\leq 
\varepsilon}\mu(y,e^{i\theta}e).$$
\end{lem}
\begin{proof}
Puisque $1$ est isol\'e (ou n'est pas) dans $U_{x,e}$, soit $0<\varepsilon<\pi$ tel que $U_{x,e}\cap\mathcal 
A_\varepsilon=\emptyset$. Dans une alg\`ebre de Jordan Banach, l'ensemble des \'el\'ements 
inversibles est ouvert, donc il existe un voisinage $\mathcal V$ de $x$ tel que 
$$\forall y\in\mathcal V\ \ y-e^{\pm i\varepsilon}e\ \text{est inversible}.$$
Alors si $y$ est une unit\'e dans $\mathcal V$, $\sigma_\varepsilon=\mathcal 
A_\varepsilon\cap U_{y,e}$ est un sous-ensemble spectral et on peut donc d\'efinir 
$q(y,e,\sigma_\varepsilon)=\mathcal G_{y,e}^{-1}(\chi_{\sigma_\varepsilon})$. Alors (cf. par 
exemple
\cite[IX, lemma 13]{DNS2})
$$p(x,e)=\int_{\nm{\lambda-1}=\varepsilon}{\frac1{2i\pi}(\lambda e-x)^{-1}d\lambda},$$
et
$$q(y,e,\sigma_\varepsilon)=\int_{\nm{\lambda-1}=\varepsilon}{\frac1{2i\pi}(\lambda 
e-y)^{-1}d\lambda},$$
ce qui montre, l'inversion \'etant continue, que si $y$ est suffisamment proche de $x$, alors
$q(y,e,\sigma_\varepsilon)$ l'est suffisamment de $p(x,e)$. 
Or si $p$ est un idempotent d'une $JB$-alg\`ebre $A$, tout idempotent $q$ dans un voisinage 
de $p$ peut s'\'ecrire 
$$q=\exp{k_v}(p)$$
o\`u $v\in A_\frac12(p)$ et $\exp {k_v}$ est un automorphisme 
de A (cf. \cite{CHIS}). Quitte \`a restreindre $\mathcal V$, on a donc, en faisant $p=p(x,e)$ et 
$q=q(y,e,\sigma_\varepsilon)$ : pour tout $y$ dans $\mathcal V$, $A_q(e)$ est isomorphe \`a 
$A_p(e)$. En particulier, si $A_p(e)$ est de rang fini alors $A_q(e)$ aussi et les rangs sont \'
egaux. 
De plus, dans ce cas, d'après le lemme précédent, l'ensemble $\sigma_\eh$ est fini. 
Supposons $\sigma_\eh=\{e^{i\theta_1},\dots,e^{i\theta_l}\}$ 
et soit 
$q_j=\mathcal{G}_{y,e}^{-1}(e^{i\theta_j})$, alors en 
faisant le calcul dans $C(U_{y,e})$, on 
voit que les $q_j$ sont des idempotents deux à deux 
orthogonaux tels que
$$q=q_1+\dots+q_l,$$
et donc $\rk(q)=\rk(q_1)+\dots+\rk(q_l).$    
\end{proof}

\section{L'indice de Maslov}

On considère dans cette partie un chemin continu $x:~\left[0,1\right]\rightarrow\Sigma$ tel que pour tout $t\in\left[0,1\right]$, $(x(t),e)$ soit une paire de Fredholm. En particulier, pour tout $t$, $1$ est une valeur propre isolée de multiplicité finie de $x(t)$ par rapport à $e$. Les résultats suivant généralisent ceux de \cite{FU} et grâce à la partie précédente les démonstrations sont semblables et parfois laissées au lecteur.
\begin{lem} 
Il existe une subdivision $t_0=0<t_1<\dots<t_N=1$ et des réels $\varepsilon_j\in\left]0,\pi\right[$, $j=1\dots N$ tels que $\forall t\in \left[t_{j-1},t_j\right]$ 
\begin{equation*}
\mu(x(t),e,\pm\varepsilon_j)=0,
\end{equation*}
et
\begin{align*}
Sp(x(t),e)\cap\mathcal{A}_{\varepsilon_j}\ \ &\text{consiste en un nombre fini de valeurs}\\
&\text{ propres isolées de multiplicités finies}.
\end{align*}
\end{lem}
\demo
La démonstration se copie sur celle de \cite[lemma 3.1]{FU} en utilisant le lemme~\ref{pert}. On l'applique à chaque $(x(t),e)$ pour $t\in\left[0,1\right]$ et on obtient des voisinages ${\mathcal{V}}_t$ et des $\varepsilon_t$. On extrait un recouvrement fini de $\left[0,1\right]$ : $$\left[0=s_0,s_0+\delta_0^{+}\right[,\dots,\left]s_i-\delta_i^{-},s_i+\delta_i^{+}\right[,\dots,\left]s_{N-1}-\delta_{N-1}^{-},s_{N-1}=1\right]$$ et on pose $$t_0=s_0=0,t_1=s_0+\delta_0^{+},\dots,t_{N-1}=s_{N-2}+\delta_{N-2}^{+},t_N=s_{N-1}=1$$ 
et 
$$\varepsilon_j=\varepsilon_{s_{j-1}}.$$
\fdemo
On dira qu'une telle subdivision $t_0=0<t_1<\dots<t_N=1$ est admissible pour les $\varepsilon_j$, $j=1,\dots,N$. Posons 
$$k(t,\varepsilon_j)=\sum_{0\leq\theta\leq\varepsilon_j}\mu(x(t),e,\theta)\quad\text{pour}\ t_{j-1}\leq t\leq t_j.$$
\begin{lem}
Soit $t_0=0<t_1<\dots<t_N=1$ une subdivision admissible pour les $\varepsilon_j$, $j=1,\dots,N$ et les $\widetilde{\varepsilon_j}$, $j=1,\dots,N$. Alors pour tout $1\leq j\leq N$,
$$k(t_j,\varepsilon_j)-k(t_{j-1},\varepsilon_j)=k(t_j,\widetilde{\varepsilon}_j)-k(t_{j-1},\widetilde{\varepsilon}_j)$$
\end{lem}
\demo
Supposons que $\varepsilon_j\geq\widetilde{\varepsilon}_j$. Alors $$k(t,\varepsilon_j)-k(t,\widetilde{\varepsilon}_j)=\sum_{\widetilde{\varepsilon}_j\leq \theta\leq \varepsilon_j}\mu(x(t),e,\theta).$$
Mais si $\gamma$ est le cercle de diamètre $\left[e^{i\varepsilon_j},e^{i\widetilde{\varepsilon}_j}\right]$, et $p_t=\frac1{2i\pi}\int_\gamma(\lambda e-x(t))^{-1}\mathrm{d}\lambda$ alors $\sum_{\widetilde{\varepsilon}_j\leq \theta\leq \varepsilon_j}\mu(x(t),e^{i\theta}e)=\rang{(p_t)}$. Mais $t\mapsto p_t$ est continue donc le rang de $p_t$ est constant.
\fdemo
\begin{pd}
La quantité 
$$ Mas({x(t)},e)=\sum_{j=1}^{N}(k(t_j,\varepsilon_j)-k(t_{j-1},\varepsilon_j)) $$ ne dépend ni des $t_j$, ni des $\varepsilon_j$, pourvu que la subdivision $t_0,\dots,t_N$ soit admissible pour les $\varepsilon_j$. On l'appelle l'indice de Maslov du chemin ${x(t)}$ par rapport au point $e$.
\end{pd}
\demo
La démonstration utilise le lemme précédent comme dans \cite[Proposition 3.3]{FU}.
\fdemo
\begin{theo}
L'indice de Maslov (par rapport à un point fixé) est additif pour la concaténation des chemins et invariant par homotopie.
\end{theo}
Le Lemme~\ref{pert} permet encore une fois d'adapter la démonstration de \cite[Theorem 3.6]{FU}.

\section{la dimension finie}

Dans cette partie on suppose $E$ de dimension finie.
Soient $x$ et $e$ dans $\Sigma$. Il existe d'uniques nombres complexes $u_1,\dots,u_k$, de module $1$ et tous distincts, et un unique système complet d'idempotents orthogonaux $c_1,\dots,c_k$ de l'algèbre de Jordan $A(e)$ tels que (cf. \cite[Proposition X.2.3 et Theorem III.1.1]{FAKO})
$$x=u_1c_1+\dots+u_kc_k.$$
L'indice de transversalité est simplement 
$$\mu(x,e,\theta)=\mu(x,e^{i\theta}e)=\sum_{j\ \text{tq}\ u_j=e^{i\theta}}\rang{(c_j)}.$$

\begin{ex}
On calcule l'indice de Maslov dans le cas du cercle pour les chemins suivants :

(i) $x(t)=e^{it\varphi}e$, où $\varphi\in [0,\pi[$.\\
On choisit $\varphi<\varepsilon<\pi$ et alors $$\mu(x(t),e^{\pm i\varepsilon}e)=0\ \text{pour}\  t\in\left[0,1\right]$$ et donc $$Mas({x(t)},e)=\sum_{0\leq\theta\leq\varepsilon}\mu(e^{i\varphi}e,e^{i\theta}e)-\sum_{0\leq\theta\leq\varepsilon}\mu(e,e^{i\theta}e)=1-1=0$$

(ii) $x(t)=e^{i(t\varphi+\psi)}e$, où $\psi\in ]0,2\pi[$ et $0<\varphi<2\pi-\psi$.\\
On choisit $0<\varepsilon<\mathrm{min}\{\psi,2\pi-(\varphi+\psi)\}$ et alors $$\mu(x(t),e^{\pm i\varepsilon}e)=0\ \text{pour}\  t\in\left[0,1\right]$$ et donc $$Mas({x(t)},e)=\sum_{0\leq\theta\leq\varepsilon}\mu(e^{i(\varphi+\psi)}e,e^{i\theta}e)-\sum_{0\leq\theta\leq\varepsilon}\mu(e^{i\psi}e,e^{i\theta}e)=0-0=0$$

(iii) $x(t)=e^{-it\varphi}e$, où $\varphi\in [0,\pi[$.\\
On choisit $\varphi<\varepsilon<\pi$ et alors $$\mu(x(t),e^{\pm i\varepsilon}e)=0\ \text{pour}\  t\in\left[0,1\right]$$ et donc $$Mas({x(t)},e)=\sum_{0\leq\theta\leq\varepsilon}\mu(e^{i\varphi}e,e^{i\theta}e)-\sum_{0\leq\theta\leq\varepsilon}\mu(e,e^{i\theta}e)=0-1=-1.$$

\end{ex}

On peut alors construire un indice pour les couples de points dans le revêtement universel $\widetilde \Sigma$ de $\Sigma$ (dont une construction se trouve dans \cite{CLKO}). Si $\widetilde\sigma$ et $\widetilde\tau$ ont pour projections respectives $\sigma$ et $\tau$ alors on pose $Mas(\widetilde\sigma,\widetilde\tau,e)=Mas({x(t)},e)$ où $x$ est n'importe quel chemin d'extrémités $\sigma$ et $\tau$ dont le relèvement d'origine $\widetilde\sigma$ se termine en $\widetilde\tau$. On note $m(\widetilde\sigma,\widetilde\tau)$ l'indice de Souriau généralisé\footnote{Dans le cas de la Lagrangienne, cet indice est en fait le double de l'indice de Souriau.} construit dans \cite{CLKO} et $\iota(\sigma_1,\sigma_2,\sigma_3)$ l'indice triple de \cite{CL}.  

On suppose désormais que $E$ est simple, autrement dit que $A(e)$ est simple (ie. ne contient pas d'idéal non trivial). Alors la composante connexe $K^{(e)}$ du groupe des automorphismes de $A(e)$ agit de manière transitive sur l'ensembles des repères de Jordan de $A(e)$ (systèmes complets d'idempotents primitifs, cf. \cite{FAKO}).

\begin{theo}
Soient $\widetilde\sigma$ et $\widetilde\tau$ dans $\widetilde \Sigma$, de projections respectives $\sigma$ et $\tau$, et soit $e$ dans $\Sigma$. Alors
\begin{equation}\tag{E}\label{E}
Mas(\widetilde\sigma,\widetilde\tau,e)=\frac12(m(\widetilde\sigma,\widetilde\tau)+\iota(e,\tau,\sigma)+\mu(\tau,e)-\mu(\sigma,e)).
\end{equation}
\end{theo}
\demo
Soient 
$$\widetilde\sigma=(\sigma=\sum e^{i\varphi_j}c_j,r\varphi)\quad \text{et}\quad 
\widetilde\tau=(\tau=\sum e^{i\phi_j}d_j,r\phi)$$ 
deux points de 
$\widetilde \Sigma$, $(c_j)$ et $(d_j)$ étant deux repères de Jordan de 
$A(e)$. Posons $\widetilde\tau'=(\tau'=\sum e^{i\phi_j}c_j,r\phi)$. Il 
existe $k\in K^{(e)}$ tel que $kd_j=c_j$, $j=1\dots r$, et soit 
$t\mapsto k_t$ un chemin dans $K^{(e)}$ tel que $k_0=id$ et $k_1=k$. 
Alors $t\mapsto (\sum e^{i\phi_j}k_td_j,r\phi)$ est un chemin dans 
$\widetilde \Sigma$ et on note $t\mapsto x(t)$ sa projection. Alors 
$Mas(\{x(t)\},e)$ est nul car $\mu(x(t),e)$ est constant. Soit $\widetilde e$ un point de $\widetilde \Sigma$ au dessus de $e$. Alors d'apr\`es la formule de Leray (cf. \cite{CLKO} 
$$m(\widetilde\tau,\widetilde\tau')+\iota(e,\tau',\tau)+\mu(\tau',e)-\mu(\tau
,e)=m(\widetilde e,\widetilde\tau')-m(\widetilde e,\widetilde\tau).$$
Comme le second membre de \eqref{E} est 
aussi additif pour la concaténation des chemins, il suffit de démontrer 
\eqref{E} en remplaçant $\widetilde\tau$ par $\widetilde\tau'$. 
Supposons que $\varphi_j\in [0,2\pi[$ et que $r\varphi=\sum\varphi_j$, 
ce qui est possible sans perte de généralité. On considère le chemin
$$t\mapsto\sum e^{i(1-t)\varphi_j}c_j$$
Son relevé d'origine $\widetilde\sigma$ se termine en $(e,0)$, et son indice de Maslov est nul puisqu'il se décompose en une succession de chemins unidimensionnels d'indices de Maslov nuls. D'autre part, si on pose $l=\#\{j\mid \varphi_j=0\ \left[2\pi\right]\}$ alors
$$m(\widetilde\sigma,(e,0))+\iota(e,e,\sigma)+\mu(e,e)-\mu(\sigma,e)=-(r-l)+0+r-l=0,$$
et donc il reste à montrer que le formule est vrai si $\widetilde\sigma=(e,0)$. Supposons que  
$\phi_j\in [0,2\pi[ $ et que $r\phi=\sum\phi_j+2k\pi$. L'indice de Maslov du chemin
$$t\mapsto\sum{e^{it\phi_j}}$$ est nul, et si $l=\#\{j\mid \phi_j=0 \left[2\pi\right]\}$ alors $$m((e,0),(\tau,\sum\phi_j))+\iota(e,\tau,e)+\mu(\tau,e)-\mu(e,e)=r-l+0+l-r=0.$$ Considérons enfin le chemin
$$t\mapsto e^{i\phi_1+2kt\pi}+\sum_{j=2}^{r}e^{i\phi_j},$$
dont le relevé d'origine $(\tau,\sum\phi_j)$ se termine en $\widetilde\tau$. Son indice de Maslov vaut $k$, tout comme le membre de droite de \eqref{E}.
\fdemo
\begin{re}
En faisant $\sigma=\tau$ dans l'équation \eqref{E} on voit que l'indice d'un lacet ne dépend pas du point par rapport auquel on le calcule.
\end{re}

\begin{re}
En faisant $\sigma=e$, on obtient $$Mas(\widetilde\sigma,\widetilde\tau,\sigma)=\frac12(m(\widetilde\sigma,\widetilde\tau)+\mu(\tau,\sigma)-r)$$
et donc on peut retrouver l'indice de Souriau, puis l'indice triple par la formule de Leray, grâce à ce nouvel indice. 
\end{re}

Un indice triple a été construit en dimension infinie par Neeb et {\O}rsted (cf. \cite{NEOR}), mais il est à valeur dans le groupe fondamental du groupe structural de $E$.
\begin{pb}
Associer une quantité numérique à cet indice, et pouvoir retrouver cet indice triple numérique grâce à l'indice pour les chemins.  
\end{pb}
\bibliographystyle{amsalpha}
\bibliography{mesref2}
\end{document}